\documentclass[10pt]{article}
\usepackage{amsfonts}
\usepackage{hyperref}
\usepackage{mathrsfs,enumerate,amssymb,amsmath,color,lineno}
\usepackage{indentfirst}
\usepackage{amsmath}
\usepackage{amssymb}
\usepackage[amsmath,thmmarks]{ntheorem}
\usepackage{graphicx}
\usepackage{tikz}

\newtheorem{definition}{\bf Definition}[section]
\newtheorem{lemma}{\bf Lemma}[section]
\newtheorem{theorem}{\bf Theorem}[section]
\newtheorem{remark}{\bf Remark}[section]

\newtheorem{example}{\bf Example}[section]

\newtheorem{proposition}{\bf Proposition}[section]

\begin{document}
\setcounter{page}{1}

\title{{\textbf{Note on additive generator pairs of overlap and grouping functions}}\thanks {Supported by
the National Natural Science Foundation of China (No.12471440)}}
\author{Li-zhi Liang\footnote{\emph{E-mail address}: 1902428258@qq.com}, Xue-ping Wang\footnote{Corresponding author. xpwang1@hotmail.com; fax: +86-28-84761502},\\
\emph{School of Mathematical Sciences, Sichuan Normal University,}\\
\emph{Chengdu 610066, Sichuan, People's Republic of China}}

\newcommand{\pp}[2]{\frac{\partial #1}{\partial #2}}
\date{}
\maketitle
\begin{quote}
{\bf Abstract} In this article, we deeply reveal the relationship between functions $\theta$ and $\vartheta$ in an overlap function additively generated by an additive generator pair ($\theta$,$\vartheta$). Then we characterize the conditions for an overlap function additively generated by the pair being a triangular norm by terms of functions $\theta$ and $\vartheta$. We also establish the conditions that an overlap function additively generated by the additive generator pair can be obtained by a distortion of a triangular norm and a (pseudo) automorphism. Finally, we dually give the related results concerned grouping functions.

\end{quote}

\section{Introduction}
Overlap and grouping functions introduced by Bustince et al. \cite{HB2010,HB2012} arise from some problems in image processing \cite{HB2007}, classification and decision making \cite{JA2010}, for which triangular norms (t-norms for short) and triangular conorms (t-conorms for short) are widely considered \cite{JF1994} but the associativity is not strongly required in the real world. As two special cases of aggregation functions \cite{GB2007,GM1986}, overlap and grouping functions are given by increasing continuous commutative bivariate functions defined over the unit square and with appropriate boundary conditions.

In recent years, overlap and grouping functions have a very rapid development both in applications and theory. In applications, overlap and grouping functions play an important role in image processing \cite{HB2010,AJ2013}, classification \cite{ME2016,ME2015} and decision making problems \cite{HB2012} when we need to assign a given object into one of two different classes which are not evidently distinct. In theory, several explorations have also been made. For example, Bedregal et al. \cite{BB2013} discussed some significant properties related to overlap and grouping functions, such as migrativity, homogeneity, idempotency and the existence of generators. Particularly, they obtained important results related to the action of automorphisms on overlap and grouping functions, and investigated the preservation of those properties and also the Lipschitz condition. More theoretic achievements are referred in \cite{BB2014,GP2015,GD2015,GP2014,DG2016}.

As we know, there exists a close relationship between overlap and grouping functions and some particular class of t-norms and t-conorms. Because a very significant way for building t-norms or t-conorms is the use of additive and multiplicative generators \cite{EP2000,OY2007,PV2005}, it seems natural for us to investigate whether such a building method are suitable for the overlap and grouping functions or not. Moreover, with an eye kept on application, the use of additive and multiplicative generators may simplify the choice of an appropriate overlap or grouping function for a given problem since we only need to consider one-variable functions instead of a bivariate one, reducing the computational complexity in this way \cite{AM2009}. Based on this idea, Dimuro et al. \cite{GD2014} supplied the notion of an additive generator pair for grouping functions. Later, they also suggested the notion of an additive generator pair for overlap functions and studied the overlap function obtained by the distortion of a positive continuous t-norm and a pseudo automorphism. Specific speaking, they gave the definition of an additive generator pair for overlap functions as follows.
\begin{definition}
	\emph{Let $\theta:[0,1]\rightarrow[0,\infty]$ and $\vartheta:[0,\infty]\rightarrow[0,1]$ be continuous and decreasing functions such that
	\item[(1)] $\theta(x)=\infty$ if and only if $x=0$;
	\item[(2)] $\theta(x)=0$ if and only if $x=1$;
	\item[(3)] $\vartheta(x)=1$ if and only if $x=0$;
	\item[(4)] $\vartheta(x)=0$ if and only if $x=\infty$.\\
	Then the function $O_{\theta,\vartheta}:[0,1]^{2}\rightarrow[0,1]$, defined by $O_{\theta,\vartheta}(x,y)=\vartheta(\theta(x)+\theta(y))$, is an overlap function. $(\theta,\vartheta)$ is called an additive generator pair of the overlap function $O_{\theta,\vartheta}$, and $O_{\theta,\vartheta}$ is said to be additively generated by the pair $(\theta,\vartheta)$.}	
\end{definition}
In 2018, Junsheng Qiao \cite{QJS2018} gave another way for building overlap and grouping functions by a multiplicative generator pair. The definition of a multiplicative generator pair is as follows.

 \begin{definition}\label{definition1.2}
 	\emph{Let $g,h:[0,1]\rightarrow[0,1]$ be two continuous and increasing functions such that
 	\item[(1)] $h(x)=0$ if and only if $x=0$;
 	\item[(2)] $h(x)=1$ if and only if $x=1$;
 	\item[(3)] $g(x)=0$ if and only if $x=0$;
 	\item[(4)] $g(x)=1$ if and only if $x=1$.\\
 	Then, the function $O_{g,h}:[0,1]^{2}\rightarrow[0,1]$, defined by $O_{g,h}(x,y)=g(h(x)h(y))$, is an overlap function. $(g,h)$ is called a multiplicative generator pair of the overlap function $O_{g,h}$, and $O_{g,h}$ is said to be multiplicatively generated by the pair $(g,h)$.}
 \end{definition}
In particular, Feng-qing Zhu \cite{FZ2022} supplied a positive answer to an open problem raised by Dimuro et al. \cite{GP2016} by constructing an example and presented the conditions under which an overlap function can be obtained by a $(\mathcal{F},T)$-distortion for a positive and continuous t-norm $T$ and a (pseudo) automorphism $\mathcal{F}$.

Notice that the condition of an additive generator pair proposed in \cite{GP2016} (resp. the condition of a multiplicative generator pair presented in \cite{QJS2018}) is sufficient but not necessary for the function $O_{\theta,\vartheta}$ (resp. $O_{g,h}$) being an overlap function. For instance, when $h$ and $g$ do not satisfy (2) and (4) of Definition \ref{definition1.2}, respectively, $O_{g,h}:[0,1]^{2}\rightarrow[0,1]$ can also be an overlap function (see Example 3.1 in \cite{QJS2018}). Thus some conditions of the additive (resp. multiplicative) generator pair can be omitted or replaced. Therefore, an interesting question naturally arise: do we characterize the functions $\theta$ and $\vartheta$ (resp. $g$ and $h$) and explore whether we can construct overlap and grouping functions by an additive generator pair $(\theta,\vartheta)$ (resp. $(g,h)$) with weaker conditions or not? This article will positively answer this problem.

The rest of this article is organized as follows: In Section \ref{section2}, we present some preliminary concepts and results. In Section \ref{section3}, we first introduce the concept of an additive generator pair ($\theta$,$\vartheta$) of an overlap function, and then investigate the relationship between functions $\theta$ and $\vartheta$. In Section \ref{section4}, we shall give the conditions for an overlap function additively generated by the pair($\theta$,$\vartheta$) being a t-norm. In Section \ref{section5}, we explore the conditions that an overlap function additively generated by the pair($\theta$,$\vartheta$) can be obtained by a distortion of a t-norm and a (pseudo) automorphism. Section \ref{section6} is devoted to introducing the concept of an additive generator pair of grouping functions and some related results. A conclusion is drawn in Section \ref{section7}.

\section{Preliminaries}\label{section2}

In this section, we recall some fundamental concepts and results that will be used in the sequel.

\begin{definition}[\cite{EP2000}]
\emph{A t-norm is a binary operator $T:[0,1]^2\rightarrow[0,1]$ such that for all $x,y,z\in[0,1]$, the following conditions are satisfied:
\begin{itemize}
	\item[(T1)] $T(x,y)=T(y,x)$;
	\item[(T2)] $T(T(x,y),z)=T(x,T(y,z))$;
	\item[(T3)] $T(x,y)\leq T(x,z)$ whenever $y\leq z$;
	\item[(T4)] $T(x,1)=x$.
\end{itemize}}
\end{definition}

A t-norm $T$ is called positive if $T(x,y)=0$ then either $x=0$ or $y=0$. A binary operator $T:[0,1]^2\rightarrow [0,1]$ is called t-subnorm if it satisfies $(T1),(T2),(T3)$ and $T(x,y)\leq \min\{x,y\}$ for all $x,y\in [0,1]$.

Let $T$ be a t-norm. Define for each $n$-tuple $(x_{1}, x_{2}, \ldots, x_{n})\in [0,1]^{n}$
$$\underset{i=1}{\overset{n}\intercal}x_{i}=T(\underset{i=1}{\overset{n-1}\intercal}x_{i}, x_{n})=T(x_{1}, x_{2}, \ldots, x_{n}).$$
If, in particular, we have $x_{1}= x_{2}= \cdots=x_{n}=x$, we shall briefly write
$$x_{T}^{(n)}=T(x,x,\ldots,x).$$

\begin{definition}[\cite{EP2000}]\label{de2.2}
\emph{Let $p,q,s,t\in[-\infty,\infty]$ with $p<q,s<t$ and $f:[p,q]\rightarrow[s,t]$ be an increasing (resp. a decreasing) function. Then the function $f^{(-1)}:[s,t]\rightarrow[p,q]$ defined by
\begin{equation*}
f^{(-1)}(y)=\sup\{x\in [p,q]\mid f(x)<y\}\,(\mbox{resp. }f^{(-1)}(y)=\sup\{x\in [p,q]\mid f(x)>y\})
\end{equation*}
is called a pseudo-inverse of the increasing (resp. decreasing) function $f$.}
\end{definition}

\begin{definition}[\cite{HB2010}]\label{definition2.3}
\emph{A bivariate function $O:[0,1]^{2}\rightarrow[0,1]$ is said to be an overlap function if it satisfies the following conditions:
\begin{itemize}
	\item[(O1)] $O$ is commutative;
	\item[(O2)] $O(x,y)=0$ if and only if $xy=0$;
	\item[(O3)] $O(x,y)=1$ if and only if $xy=1$;
	\item[(O4)] $O$ is increasing;
	\item[(O5)] $O$ is continuous.
\end{itemize}}
\end{definition}

\begin{definition}[\cite{HB2012}]\label{definition2.4}
	\emph{A bivariate function $G:[0,1]^{2}\rightarrow[0,1]$ is said to be a grouping function if it satisfies the following conditions:
	\begin{itemize}
		\item[(G1)] $G$ is commutative;
		\item[(G2)] $G(x,y)=0$ if and only if $x=y=0$;
		\item[(G3)] $G(x,y)=1$ if and only if $x=1$ or $y=1$;
		\item[(G4)] $G$ is increasing;
		\item[(G5)] $G$ is continuous.
	\end{itemize}}
\end{definition}

\begin{definition}[\cite{GP2016}]\label{definition2.5}
	\emph{A function $\mathcal{F}:[0,1]\rightarrow[0,1]$ is said to be a pseudo automorphism if the following conditions hold:
	\begin{itemize}
		\item[($\mathcal{F}$1)] $\mathcal{F}$ is increasing;
		\item[($\mathcal{F}$2)] $\mathcal{F}$ is continuous;
		\item[($\mathcal{F}$3)] $\mathcal{F}(x)=1$ if and only if $x=1$;
		\item[($\mathcal{F}$4)] $\mathcal{F}(x)=0$ if and only if $x=0$.
	\end{itemize}}
\end{definition}
An automorphism $\mathcal{F}:[0,1]\rightarrow[0,1]$ is a strictly increasing pseudo automorphism.

\begin{proposition}[\cite{GP2016}]\label{prop1}
	Let $\mathcal{F}:[0,1]\rightarrow[0,1]$ be a pseudo automorphism. Then for every positive and continuous t-norm $T:[0,1]^{2}\rightarrow[0,1]$, the function $O_{\mathcal{F},T}:[0,1]^{2}\rightarrow[0,1]$, given by
	\begin{equation*}
		O_{\mathcal{F},T}(x,y)=\mathcal{F}(T(x,y)),
	\end{equation*}
    is an overlap function.
     \end{proposition}

     The function $O_{\mathcal{F},T}$ in Proposition \ref{prop1} is called an overlap function obtained by the distortion of the t-norm $T$ by the pseudo automorphism $\mathcal{F}$ or an overlap function obtained by the $(\mathcal{F},T)$-distortion.
\section{Additive generator pairs of overlap functions}\label{section3}

This section first introduces the concept of an additive generator pair ($\theta$,$\vartheta$) of an overlap function, and then it deeply reveals the relationship between functions $\theta$ and $\vartheta$ in the overlap function additively generated by the pair.

Similar to Definition 3.1 of \cite{QJS2018}, we may give the following one.
\begin{definition}\label{definition1}
	\emph{Let $\theta:[0,1]\rightarrow[0,\infty]$ and $\vartheta:[0,\infty]\rightarrow[0,1]$ be continuous and decreasing functions, respectively. If the bivariate function $O_{\theta,\vartheta}:[0,1]^{2}\rightarrow[0,1]$, given by
	$$O_{\theta,\vartheta}(x,y)=\vartheta(\theta(x)+\theta(y)),$$
	is an overlap function, then $(\theta,\vartheta)$ is called an additive generator pair of the overlap function $O_{\theta,\vartheta}$ and $O_{\theta,\vartheta}$ is said to be additively generated by the pair $(\theta,\vartheta)$.}
\end{definition}

Let $\theta:[0,1]\rightarrow[0,\infty]$ be a function. For any $a\in [0, \infty)$, define functions $2\theta:[0,1]\rightarrow[0,\infty]$ and $\theta+\frac{a}{2}:[0,1]\rightarrow[0,\infty]$, respectively, by $$(2\theta)(x)=2\theta(x)$$ for any $x\in[0,1]$ and $$(\theta+\frac{a}{2})(x)=\theta(x)+\frac{a}{2}$$ for any $x\in[0,1]$. Then the following theorem characterizes the functions $\theta:[0,1]\rightarrow[0,\infty]$ and $\vartheta:[0,\infty]\rightarrow[0,1]$, respectively, whenever $O_{\theta,\vartheta}$ is an overlap function additively generated by the pair ($\theta$,$\vartheta$).
\begin{theorem}\label{theorem3.2}
	Let $\theta:[0,1]\rightarrow[0,\infty]$ and $\vartheta:[0,\infty]\rightarrow[0,1]$ be continuous and decreasing functions, respectively, and $O_{\theta,\vartheta}:[0,1]^2\rightarrow[0,1]$ be an overlap function additively generated by the pair $(\theta,\vartheta)$. Then the following two statements hold:
	\begin{itemize}
		\item[(i)] $\theta(x)=\infty$ if and only if $x=0$;
		\item[(ii)] $\vartheta(x)=0$ if and only if $x=\infty$.
	\end{itemize}
\end{theorem}
\begin{proof}
	(i) ($\Leftarrow$) If $x=0$, then we affirm that $\theta(x)=\infty$. Indeed, if $\theta(0)<\infty$ then by Definition \ref{definition2.3}, $\vartheta(u)\le\vartheta(2\theta(0))=\vartheta(\theta(0)+\theta(0))=O_{\theta,\vartheta}(0,0)=0$, i.e.,
\begin{equation}\label{eq1}
\vartheta(u)=0
\end{equation}
 for any $u\in[2\theta(0),\infty]$.	Note that $2\theta:[0,1]\rightarrow[0,\infty]$ is continuous since $\theta:[0,1]\rightarrow[0,\infty]$ is continuous. In the following, we prove that $\vartheta(y)\ne0$ for any $y\in[0,2\theta(0))$ by distinguishing two cases.
	
	Case 1. If $y\in(2\theta(1),2\theta(0))$, then there exists a $z\in(0,1)$ such that $y=2\theta(z)$ by the continuity of $2\theta$. Assume that $\vartheta(y)=0$. Then one has that
	\begin{align*}
		O_{\theta,\vartheta}(z,z)=&\vartheta(\theta(z)+\theta(z))\\
		=&\vartheta(2\theta(z))\\
		=&\vartheta(y)\\
		=&0,
	\end{align*}
	which contradicts (O2) of Definition\ref{definition2.3}. Therefore,
	\begin{equation}\label{equation(1)}
		\vartheta(y)\ne0
	\end{equation}
	for any $y\in(2\theta(1),2\theta(0))$.

	Case 2. If $y\in[0,2\theta(1)]$ and $\vartheta(y)=0$, then it is obvious that $\vartheta(u)=0$ for any $u\in(2\theta(1),2\theta(0))$ since $\vartheta$ is a decreasing function, contrary to \eqref{equation(1)}. Therefore, $\vartheta(y)\neq 0$ for any $y\in[0,2\theta(1)]$.

    Cases 1 and 2 imply that $\vartheta(y)\ne0$ for any $y\in[0,2\theta(0))$, which together with \eqref{eq1} implies that
    \begin{equation}\label{eq2}\vartheta(u)=0 \mbox{ if and only if }u\in[2\theta(0),\infty].
    \end{equation}
	On the other hand, since $O_{\theta,\vartheta}$ is an overlap function, by (O2) of Definition \ref{definition2.3} we have that
	$$O_{\theta,\vartheta}(0,y)=\vartheta(\theta(0)+\theta(y))=0$$
 for any $y\in[0,1]$.
	Thus from \eqref{eq2} we get $\theta(0)+\theta(y)\ge2\theta(0)$, i.e., $\theta(y)\ge\theta(0)$ for any $y\in[0,1]$. In particular, $\theta(0)\le\theta(1)$. Thus $\theta(1)=\theta(0)$ since $\theta(1)\le\theta(0)$. Furthermore, we conclude that
	\begin{align*}
		O_{\theta,\vartheta}(1,1)=&\vartheta(\theta(1)+\theta(1))\\
		=&\vartheta(\theta(0)+\theta(0))\\
		=&O_{\theta,\vartheta}(0,0)\\
		=&0\quad\mbox{ by (O2) of Definition \ref{definition2.3}}
	\end{align*}
	which contradicts (O3) of Definition \ref{definition2.3}. Therefore, $\theta(0)=\infty$.
	
	($\Rightarrow$) If $\theta(x)=\infty$ then we assert that $x=0$. In fact, if there exists a $u\in(0,1]$ such that $\theta(u)=\infty$, then
	\begin{align*}
		O_{\theta,\vartheta}(u,u)=&\vartheta(\theta(u)+\theta(u))\\
		=&\vartheta(\infty),
	\end{align*}
Since $\theta(0)=\infty$, by (O2) of Definition \ref{definition2.3}, we further have
	\begin{equation}\label{equation(2)}
		O_{\theta,\vartheta}(u,u)=\vartheta(\infty)=\vartheta(\theta(0)+\theta(0)))=O_{\theta,\vartheta}(0,0)=0,
	\end{equation}
		contrary to (O2) of Definition \ref{definition2.3}.
	
	(ii) ($\Leftarrow$) If $x=\infty$, then by \eqref{equation(2)}, we know that $\vartheta(x)=0$.
	
	($\Rightarrow$) From (i), we know that $\theta(x)=\infty$ if and only if $x=0$. Thus if $\vartheta(x)=0$ then we assert that $x=\infty$. Otherwise, there are two cases as follows.

Case a. If there exists a $u\in(2\theta(1),2\theta(0))$ such that $\vartheta(u)=0$, then there exists a $z\in(0,1)$ such that $u=2\theta(z)$ by the continuity of $2\theta$. Thus
	\begin{align*}
		O_{\theta,\vartheta}(z,z)=&\vartheta(\theta(z)+\theta(z))\\
		=&\vartheta(u)\\
		=&0,
	\end{align*}
	which is contrary to (O2) of Definition \ref{definition2.3}. Consequently,
	\begin{equation}\label{equation(3)}
		\vartheta(y)\ne0
	\end{equation}
for any $y\in(2\theta(1),2\theta(0))$.
	
Case b. If there exists a $u\in[0,2\theta(1)]$ such that $\vartheta(u)=0$, then $\vartheta(y)\le\vartheta(u)=0$ for any $y\in(2\theta(1),\infty)$, contrary to \eqref{equation(3)}.

Therefore, from Cases a and b, $\vartheta(u)\neq 0$ for any $u\in[0, \infty)$, i.e., $x=\infty$.
\end{proof}

The following theorem shows the conditions that a pair $(\theta,\vartheta)$ can additively generate overlap functions.
\begin{theorem}\label{theorem3.1}
	For a given $a\in[0,\infty)$, let $\theta:[0,1]\rightarrow[0,\infty]$ and $\vartheta:[0,\infty]\rightarrow[0,1]$ be continuous and decreasing functions such that
	\begin{itemize}
		\item[(1)] $\theta(x)=\infty$ if and only if $x=0$;
		\item[(2)] $\theta(x)=\frac{a}{2}$ if and only if $x=1$;
		\item[(3)] $\vartheta(x)=1$ if and only if $x\in[0,a]$;
		\item[(4)] $\vartheta(x)=0$ if and only if $x=\infty$.
	\end{itemize}
	Then, the function $O_{\theta,\vartheta}:[0,1]^{2}\rightarrow[0,1]$, defined by
	\begin{equation*}
		O_{\theta,\vartheta}(x,y)=\vartheta(\theta(x)+\theta(y)),
	\end{equation*}
	is an overlap function.
\end{theorem}
\begin{proof}
	We only need to show that Definition \ref{definition2.3} holds. Indeed, the conditions (O1) and (O5) are obvious. So that we just verify the other conditions as follows.

	(O2) First note that from (4), $O_{\theta,\vartheta}(x,y)=0\ \text{if and only if}\ \vartheta(\theta(x)+\theta(y))=0\ \text{if and only if}\
		\theta(x)+\theta(y)=\infty$.

Now, suppose that $O_{\theta,\vartheta}(x,y)=0$. We claim that $xy=0$. In fact, if $xy\ne0$. Then $x\ne0$ and $y\ne0$. Thus from (1), we have $\theta(x)+\theta(y)\ne\infty$, a contradiction.

 Conversely, if $xy=0$, then $x=0$ or $y=0$, say $x=0$, thus from (1) we have $\theta(x)=\infty$. This follows that $\theta(x)+\theta(y)=\infty$, which implies that $O_{\theta,\vartheta}(x,y)=0$.

    (O3) First note that from (3), $O_{\theta,\vartheta}(x,y)=1\ \text{if and only if}\ \vartheta(\theta(x)+\theta(y))=1\ \text{if and only if}\ \theta(x)+\theta(y)\in[0,a]$.
		
Now, suppose that $O_{\theta,\vartheta}(x,y)=1$. Then we assert that $xy=1$. In fact, if $xy\ne1$. Then $x\ne1$ or $y\ne1$, say $x\ne1$, thus from (2) we have $\theta(x)>\theta(1)=\frac{a}{2}$. Then $\theta(x)+\theta(y)>\frac{a}{2}+\theta(y)\ge\frac{a}{2}+\frac{a}{2}=a$, which is in conflict  with $\theta(x)+\theta(y)\in[0,a]$.

Conversely, if $xy=1$, then $x=1$ and $y=1$. Thus from (2) we have that $\theta(x)+\theta(y)=\theta(1)+\theta(1)=\frac{a}{2}+\frac{a}{2}=a$, i.e., $\theta(x)+\theta(y)\in[0,a]$. Therefore, $O_{\theta,\vartheta}(x,y)=1$.

    (O4) For any $y,z\in[0,1]$ with $y\le z$, $\theta(y)\ge\theta(z)$ since $\theta$ is decreasing, which together with $\vartheta$ being decreasing means that
\begin{equation*}
	O_{\theta,\vartheta}(x,y)=\vartheta(\theta(x)+\theta(y))\le\vartheta(\theta(x)+\theta(z))=O_{\theta,\vartheta}(x,z),
\end{equation*}
i.e., (O4) holds.
\end{proof}

The following two examples illustrate Theorem \ref{theorem3.1}.
\begin{example}\label{example3.1}
	\emph{For a given $a\in[0,\infty)$, consider the functions $\theta:[0,1]\rightarrow[0,\infty]$ and $\vartheta:[0,\infty]\rightarrow[0,1]$ defined by
	$$\theta(x)=\begin{cases}
			\frac{a}{2}-lnx & \mbox{if}\ x\neq0,\\
		\infty & \mbox{if}\ x=0
			\end{cases}$$
	and
	$$\vartheta(x)=\begin{cases}
			1        &\mbox{if}\ x\in[0,a),\\
		e^{-x+a} &\mbox{if}\ x\in[a,\infty],
	\end{cases}$$
respectively. Then it is easy to see that both $\theta$ and $\vartheta$ are continuous and decreasing functions satisfying the conditions of Theorem \ref{theorem3.1}. Therefore, $O_{\theta,\vartheta}(x,y)=xy$ is an overlap function. In particular, $O_{\theta,\vartheta}=T_p$.}
\end{example}

\begin{example}\label{example3.2}
	\emph{For a given $a\in[0,\infty)$, consider the functions $\theta:[0,1]\rightarrow[0,\infty]$ and $\vartheta:[0,\infty]\rightarrow[0,1]$, defined, respectively, by
	$$\theta(x)=
	\begin{cases}
		\frac{a}{2}-lnx & \mbox{if}\ x\neq0,\\
		\infty & \mbox{if}\ x=0
	\end{cases}	
		$$
	and
	$$\vartheta(x)=
	\begin{cases}
		1& \mbox{if}\ x\in[0,a),\\
		\frac{a}{x}& \mbox{if}\ x\in[a,\infty].
	\end{cases}	$$
	Then we can see that both $\theta$ and $\vartheta$ are continuous and decreasing functions, satisfying the conditions of Theorem \ref{theorem3.1}.}
	
	\emph{In a simple computation, we have \begin{equation*}
		O_{\theta,\vartheta}(x,y)=\left\{
		\begin{aligned}
			&0&\mbox{otherwise},\\
			&\frac{a}{a-\ln xy}&\mbox{if}\ xy\ne0.
		\end{aligned}
		\right.
	\end{equation*}}
\end{example}
Therefore, by Theorem \ref{theorem3.1}, $O_{\theta,\vartheta}$ is an overlap function that is not associative, and 1 is not its neutral element.

The following two propositions reveal that the conditions (2) and (3) of Theorem \ref{theorem3.1} are closely associated whenever $O_{\theta,\vartheta}:[0,1]^{2}\rightarrow[0,1]$ is an overlap function additively generated by the pair $(\theta,\vartheta)$.
\begin{proposition}\label{proposition3.1}
	For a given $a\in[0,\infty)$, let $\theta:[0,1]\rightarrow[0,\infty]$ and $\vartheta:[0,\infty]\rightarrow[0,1]$ be two continuous and decreasing functions satisfying
	\begin{itemize}
		\item[(i)] $x\in[0,a]$ if and only if $\vartheta(x)=1$, and
		\item[(ii)] a binary function $O_{\theta,\vartheta}:[0,1]^{2}\rightarrow[0,1]$ with $O_{\theta,\vartheta}(x,y)=\vartheta(\theta(x)+\theta(y))$ is an overlap function.
	\end{itemize}
	Then, $\theta(x)=\frac{a}{2}$ if and only if $x=1$.	
\end{proposition}
\begin{proof}
	If $\theta(x)=\frac{a}{2}$, then from (i) we have $\vartheta(\theta(x)+\theta(x))=\vartheta(a)=1$. Hence $O_{\theta,\vartheta}(x,x)=1$. Since $O_{\theta,\vartheta}$ is an overlap function, by (O3) of Definition \ref{definition2.3}, we have that $x=1$.
	
	In the converse implication, if $x=1$, then $O_{\theta,\vartheta}(1,1)=1$ since $O_{\theta,\vartheta}$ is an overlap function, i.e., $\vartheta(\theta(1)+\theta(1))=1$. Hence, it follows from (i) that $\theta(1)\in[0,\frac{a}{2}]$. We assert that $\theta(1)=\frac{a}{2}$. In fact, if $\theta(1)<\frac{a}{2}$, then there exists an $x_{0}\in[0,1)$ such that $\theta(x_{0})=\frac{a}{2}$ since $\theta$ is continuous and $\theta(0)=\infty$, thus $\vartheta(\theta(x_{0})+\theta(x_{0}))=\vartheta(a)=1$, which violates (O3) of Definition \ref{definition2.3}. Therefore, $\theta(1)=\frac{a}{2}$.
\end{proof}

\begin{proposition}\label{proposition3.2}
	For a given $a\in[0,\infty)$, let $\theta:[0,1]\rightarrow[0,\infty]$ and $\vartheta:[0,\infty]\rightarrow[0,1]$ be two continuous and decreasing functions such that
	\begin{itemize}
		\item[(i)] $\theta(x)=\frac{a}{2}$ if and only if $x=1$, and
		\item[(ii)] a binary function $O_{\theta,\vartheta}:[0,1]^{2}\rightarrow[0,1]$ with $O_{\theta,\vartheta}(x,y)=\vartheta(\theta(x)+\theta(y))$ is an overlap function.
	\end{itemize}
	Then, $x\in[0,a]$ if and only if $\vartheta(x)=1$.
\end{proposition}
\begin{proof}
	Necessity. If $x\in[0,a]$, then $\vartheta(x)\ge\vartheta(a)$ since $\vartheta$ is decreasing. Meanwhile, from (i) we have $\vartheta(a)=\vartheta(\theta(1)+\theta(1))=O_{\theta,\vartheta}(1,1)=1$ since $O_{\theta,\vartheta}$ is an overlap function.	Thus $1\ge\vartheta(x)\ge\vartheta(a)=1$ for any $x\in[0,a]$, i.e., $\vartheta(x)=1$ for any $x\in[0,a]$.
	
	Sufficiency. If $\vartheta(x)=1$, Then we claim $x\in[0,a]$. Indeed, if there exists an $x_{0}\in(a,\infty]$ such that $\vartheta(x_{0})=1$, then there exist two elements $y,z\in[0,1]$ such that $\theta(y)+\theta(z)=x_0$ since $\theta$ is continuous. Thus $yz<1$ and $\vartheta(\theta(y)+\theta(z))=\vartheta(x_{0})=1$, contrary to (O3) of Definition \ref{definition2.3} since $O_{\theta,\vartheta}$ is an overlap function.
\end{proof}
\begin{remark} \emph{Compared with the additive generator pair of an overlap function in \cite{GP2016} (see Corollary 4.1 in \cite{GP2016}), the additive generator pair of an overlap function appeared in Theorem \ref{theorem3.1} gives a way to construct overlap functions with more relaxed boundary conditions. Indeed, if $a=0$, then Theorem \ref{theorem3.1} coincides with Corollary 4.1 in \cite{GP2016}.}
\end{remark}\begin{remark} \emph{From Propositions \ref{proposition3.1} and \ref{proposition3.2}, for a given $a\in[0,\infty)$, let $\theta:[0,1]\rightarrow[0,\infty]$ and $\vartheta:[0,\infty]\rightarrow[0,1]$ be two continuous and decreasing functions such that a binary function $O_{\theta,\vartheta}:[0,1]^{2}\rightarrow[0,1]$ with $O_{\theta,\vartheta}(x,y)=\vartheta(\theta(x)+\theta(y))$ is an overlap function. Then the conditions (2) and (3) of Theorem \ref{theorem3.1} are equivalent.}	
\end{remark}

\section{Conditions that an overlap function is a t-norm}\label{section4}
From Examples \ref{example3.1} and \ref{example3.2}, one easily sees that not every overlap function additively generated by the pair ($\theta$,$\vartheta$) is a t-norm. Then a natural question arises: what are the conditions for an overlap function additively generated by the pair ($\theta$,$\vartheta$) being a t-norm? This section will answer this problem.

First, we have the following proposition.
\begin{proposition}\label{proposition4.2}
	Let $\theta:[0,1]\rightarrow[0,\infty]$ and $\vartheta:[0,\infty]\rightarrow[0,1]$ be two continuous and decreasing functions, and let $O_{\theta,\vartheta}:[0,1]^{2}\rightarrow[0,1]$ be an overlap function additively generated by the pair $(\theta,\vartheta)$ with $\theta(1)=\frac{a}{2}$ where $a\in[0,\infty)$. Then the following two statements are equivalent:
	\begin{itemize}
		\item[(1)] $\vartheta\circ(\theta+\frac{a}{2})=\emph{id}_{[0,1]}$;
		\item[(2)] $1$ is a neutral element of $O_{\theta,\vartheta}$.
	\end{itemize}
\end{proposition}
\begin{proof}
	$(1)\Rightarrow(2)$ If $\vartheta\circ(\theta+\frac{a}{2})=\mbox{id}_{[0,1]}$, then
	\begin{align*}
		O_{\theta,\vartheta}(x,1)=&\vartheta(\theta(x)+\theta(1))\\
		=&\vartheta(\theta(x)+\frac{a}{2})\\
		=&x
	\end{align*} for each $x$$\in[0,1]$ since $\theta(1)=\frac{a}{2}$.
	Therefore, $1$ is a neutral element of $O_{\theta,\vartheta}$.
	
	$(2)\Rightarrow(1)$ If $1$ is a neutral element of $O_{\theta,\vartheta}$, then
	\begin{align*}
		\vartheta((\theta(x)+\frac{a}{2}))=&\vartheta((\theta(x)+\theta(1)))\\
		=&O_{\theta,\vartheta}(x,1)\\
		=&x
	\end{align*} for any $x\in[0,1]$ since $\theta(1)=\frac{a}{2}$. Therefore, $\vartheta\circ(\theta+\frac{a}{2})=\mbox{id}_{[0,1]}$.
\end{proof}

Before presenting results that are associating a t-norm with an overlap function additively generated by the pair ($\theta$,$\vartheta$), we need the following lemma.
\begin{lemma}[\cite{HB2010}]\label{lemma4.1}
	$O:[0,1]^{2}\rightarrow[0,1]$ is an associative overlap function if and only if $O$ is a continuous and positive t-norm.
\end{lemma}

\begin{proposition}\label{proposition4.1}
	For a given $a\in[0,\infty)$, let $\theta:[0,1]\rightarrow[0,\infty]$ be a continuous and strictly decreasing function such that:
	\begin{itemize}
		\item[(1)] $\theta(x)=\infty$ if and only if $x=0$;
		\item[(2)] $\theta(x)=\frac{a}{2}$ if and only if $x=1$.
	\end{itemize}
	If a function $\vartheta:[0,\infty]\rightarrow[0,1]$ satisfies $\vartheta(x)=(\theta+\frac{a}{2})^{(-1)}(x)$ for any $x\in [0,1]$, then a binary function $O_{\theta,\vartheta}:[0,1]^{2}\rightarrow[0,1]$ with $O_{\theta,\vartheta}(x,y)=\vartheta(\theta(x)+\theta(y))$ is a positive t-norm.
\end{proposition}
\begin{proof} Since $\theta$ is strictly decreasing, $\theta+\frac{a}{2}$ is also strictly decreasing. Then $\vartheta=(\theta+\frac{a}{2})^{(-1)}$ is continuous and decreasing. In the following, we verify that the $O_{\theta,\vartheta}:[0,1]^{2}\rightarrow[0,1]$ is an overlap function according to Theorem \ref{theorem3.1}. From the hypothesis it is clearly that both (1) and (2) of Theorem \ref{theorem3.1} hold. Next, we prove that (3) and (4) of Theorem \ref{theorem3.1} are true, respectively.
	
	(3) From Definition \ref{de2.2}, $\vartheta(x)=1\ \text{if and only if }\ (\theta+\frac{a}{2})^{(-1)}(x)=1\
		\text{if and only if}\ \sup\{z\in[0,1]\mid\theta(z)+\frac{a}{2}>x\}=1$.	
	
	Suppose that $\vartheta(x)=1$. If $x\in (a,\infty]$, then there exists a $u\in[0,1)$ such that $x=\theta(u)+\frac{a}{2}$ since $\theta$ is continuous and strictly decreasing. Thus $\vartheta(x)=(\theta+\frac{a}{2})^{(-1)}(x)=\sup\{z\in[0,1]\mid\theta(z)+\frac{a}{2}>\theta(u)+\frac{a}{2}\}=u\ne1$, which is contrary to $\vartheta(x)=1$. Therefore, $x\in[0,a]$.
	
	Conversely, if $x=a$, then we get
\begin{align*}
\vartheta(a)=&(\theta+\frac{a}{2})^{(-1)}(a)\\
             =&\sup\{z\in[0,1]\mid\theta(z)+\frac{a}{2}>a\}\\
             =&\sup\{z\in[0,1]\mid\theta(z)>\frac{a}{2}\}\\
             =&1.
\end{align*}
 Thus $1\ge\vartheta(x)\ge\vartheta(a)=1$, i.e., $\vartheta(x)=1$ for any $x\in[0,a]$.
	
	(4) From Definition \ref{de2.2}, $\vartheta(x)=0\ \text{if and only if}\ (\theta+\frac{a}{2})^{(-1)}(x)=0\ \text{if and only if}\
		\sup\{z\in[0,1]\mid\theta(z)+\frac{a}{2}>x\}=0$.	
	
	$(\Rightarrow)$ Supposing that $\vartheta(x)=0$, we assert that $x=\infty$. Otherwise, there exists an $x_{0}\in[0,\infty)$ such that $\vartheta(x_{0})=0$. Then, because of the continuity of $\theta$, we have
\begin{align*}
\vartheta(x_0)&=(\theta+\frac{a}{2})^{(-1)}(x_0)\\
              &=\sup\{z\in[0,1]\mid\theta(z)+\frac{a}{2}>x_{0}\}\\
              &\ne0
\end{align*}
since $\theta(x)=\infty$ if and only if $x=0$, which violates $\vartheta(x_{0})=0$.	
	
	$(\Leftarrow)$ If $x=\infty$, then $\vartheta(x)=\sup\{z\in[0,1]\mid\theta(z)+\frac{a}{2}>\infty\}=0$.
	
	Therefore, from Theorem \ref{theorem3.1}, $O_{\theta,\vartheta}$ is an overlap function.
	
Now, we prove the associativity of $O_{\theta,\vartheta}$.

For each $x,y,z\in[0,1]$, there exists an $w\in[0,1]$ such that $\theta(x)+\theta(y)=\theta(w)+\frac{a}{2}$ since $\theta$ is continuous and $\theta(x)+\theta(y)\ge\theta(x)+\frac{a}{2}$. On the other hand, $\theta+\frac{a}{2}$ is strictly decreasing since $\theta$ is strictly decreasing, which means $(\theta+\frac{a}{2})^{(-1)}(\theta(w)+\frac{a}{2})=w$. Thus
 \begin{align*}
 	 O_{\theta,\vartheta}(O_{\theta,\vartheta}(x,y),z)=&O_{\theta,\vartheta}((\theta+\frac{a}{2})^{(-1)}(\theta(x)+\theta(y)),z)\\
 =&(\theta+\frac{a}{2}) ^{(-1)}(\theta((\theta+\frac{a}{2}) ^{(-1)}(\theta(x)+\theta(y)))+\theta(z))\\
 	 =&(\theta+\frac{a}{2}) ^{(-1)}(\theta((\theta+\frac{a}{2}) ^{(-1)}(\theta(w)+\frac{a}{2}))+\theta(z))\\
 	 =&(\theta+\frac{a}{2})^{(-1)}(\theta(w)+\theta(z))\\
 	 =&(\theta+\frac{a}{2})^{(-1)}(\theta(x)+\theta(y)+\theta(z)-\frac{a}{2}).
 \end{align*}
 Analogously, $O_{\theta,\vartheta}(x,O_{\theta,\vartheta}(y,z))=(\theta+\frac{a}{2})^{(-1)}(\theta(x)+\theta(y)+\theta(z)-\frac{a}{2})$. Therefore, $O_{\theta,\vartheta}(x,O_{\theta,\vartheta}(y,z))= O_{\theta,\vartheta}(O_{\theta,\vartheta}(x,y),z)$.

Finally, by Lemma \ref{lemma4.1}, $O_{\theta,\vartheta}$ is a continuous and positive t-norm.
\end{proof}

Notice that we can directly prove Proposition \ref{proposition4.1} by the definition of the t-norm.
\begin{proposition}\label{proposition4.3}
	Let $\theta:[0,1]\rightarrow[0,\infty]$ and $\vartheta:[0,\infty]\rightarrow[0,1]$ be two continuous and decreasing functions, and $O_{\theta,\vartheta}:[0,1]^{2}\rightarrow[0,1]$ be an overlap function additively generated by the pair $(\theta,\vartheta)$ with $1$ as neutral element. If for a given $a\in [0,\infty)$, $\theta$ satisfies the condition (2) of Theorem \ref{theorem3.1}, then $O_{\theta,\vartheta}$ is associative.
\end{proposition}
\begin{proof}
	 From the condition (2) of Theorem \ref{theorem3.1} we get $\theta(x)=\frac{a}{2}$ if and only if $x=1$. Then we have
	\begin{equation}\label{eq6}
		O_{\theta,\vartheta}(x,1)=\vartheta(\theta(x)+\theta(1))=\vartheta(\theta(x)+\frac{a}{2})=x
	\end{equation}
	for any $x\in[0,1]$ since $1$ is the neutral element of $O_{\theta,\vartheta}$. Moreover, because $\theta$ is decreasing, $\theta(x)+\theta(y)\ge\theta(x)+\frac{a}{2}$ for any $x,y\in[0,1]$. Thus $\theta(x)+\theta(y)\in \mbox{Ran}(\theta+\frac{a}{2})$ since $\theta$ is continuous, i.e., there exists an $w\in[0,1]$ such that $\theta(x)+\theta(y)=\theta(w)+\frac{a}{2}$.
	It follows that for any $x,y, z\in[0,1]$,
	\begin{align*}
		O_{\theta,\vartheta}(O_{\theta,\vartheta}(x,y),z)=&\vartheta(\theta(O_{\theta,\vartheta}(x,y))+\theta(z))\\
		=&\vartheta(\theta(\vartheta(\theta(x)+\theta(y)))+\theta(z))\\
		=&\vartheta(\theta(\vartheta(\theta(w)+\frac{a}{2}))+\theta(z))\\
		=&\vartheta(\theta(w)+\theta(z))\quad\mbox{ by \eqref{eq6}}\\
		=&\vartheta(\theta(x)+\theta(y)+\theta(z)-\frac{a}{2}).
	\end{align*}
Similarly, we have $O_{\theta,\vartheta}(x,O_{\theta,\vartheta}(y,z))=\vartheta(\theta(x)+\theta(y)+\theta(z)-\frac{a}{2})$ for any $x,y, z\in[0,1]$. Therefore, $O_{\theta,\vartheta}(O_{\theta,\vartheta}(x,y),z)=O_{\theta,\vartheta}(x,O_{\theta,\vartheta}(y,z))$ for any $x,y, z\in[0,1]$.
\end{proof}

 Lemma \ref{lemma4.1}, Propositions \ref{proposition4.2} and \ref{proposition4.3} imply the following theorem.
\begin{theorem}
	Let $\theta:[0,1]\rightarrow[0,\infty]$ and $\vartheta:[0,\infty]\rightarrow[0,1]$ be two continuous and decreasing functions and, $O_{\theta,\vartheta}:[0,1]^{2}\rightarrow[0,1]$ be an overlap function additively generated by the pair $(\theta,\vartheta)$. If for a given $a\in [0,\infty)$, $\theta$ satisfies the condition (2) of Theorem \ref{theorem3.1}, then the following three statements are equivalent:
	\begin{itemize}
		\item[(1)] $O_{\theta,\vartheta}$ is a t-norm;
		\item[(2)] $1$ is a neutral element of $O_{\theta,\vartheta}$;
		\item[(3)] $\vartheta\circ(\theta+\frac{a}{2})=\emph{id}_{[0,1]}$.
	\end{itemize}
\end{theorem}

\section{Conditions that overlap functions additively generated by the pair ($\theta$,$\vartheta$) can be obtained by a distortion} \label{section5}
Dimuro et al. \cite{GP2016} introduced the definition of a (pseudo) automorphism to discuss the conditions on which an overlap function can be obtained by a distortion of a t-norm and a (pseudo) automorphism. This section investigates the conditions that an overlap function additively generated by the pair ($\theta$,$\vartheta$) can be obtained by a distortion of a t-norm and a (pseudo) automorphism.

\begin{proposition}\label{proposition5.1}
	Let $\theta:[0,1]\rightarrow[0,\infty]$ and $\vartheta:[0,\infty]\rightarrow[0,1]$ be two continuous and decreasing functions and, $O_{\theta,\vartheta}:[0,1]^{2}\rightarrow[0,1]$ be an overlap function additively generated by the pair $(\theta,\vartheta)$ that satisfies one of the following two conditions: for a given $a\in [0,\infty)$,
\begin{itemize}
		\item[(1)] $\theta(x)=\frac{a}{2}$ if and only if $x=1$;
		\item[(2)] $\vartheta(x)=1$ if and only if $x\in[0,a]$.
\end{itemize}
Then the overlap function $O_{\theta,\vartheta}$ can be obtained by a $(\mathcal{F},T_{sub})$-distortion for a (pseudo) automorphism $\mathcal{F}$ and a t-subnorm $T_{sub}$.
\end{proposition}
\begin{proof}
	Define functions $T_{sub}:[0,1]^2\rightarrow[0,1]$ and $\mathcal{F}:[0,1]\rightarrow[0,1]$, respectively, by \begin{equation}\label{eq06}T_{sub}(x,y)=(\theta+\frac{a}{2})^{(-1)}(\theta(x)+\theta(y))\end{equation} for any $x,y\in[0,1]$ and \begin{equation}\label{eq07}\mathcal{F}(x)=\vartheta(\theta(x)+\frac{a}{2})\end{equation} for any $x\in[0,1]$.
The rest of the proof is completed by three steps as below.\\
Step 1. We prove that $T_{sub}$ is a t-subnorm.
\begin{itemize}
		\item[(\romannumeral1)] The commutativity of $T_{sub}$ is obvious.
		\item[(\romannumeral2)] The associativity of $T_{sub}$. Because $\theta$ is decreasing, from the hypothesis and Proposition \ref{proposition3.1} $\theta(x)+\theta(y)\ge\theta(x)+\frac{a}{2}$ for each $x,y\in[0,1]$ since $O_{\theta,\vartheta}$ is an overlap function. Then $\theta(x)+\theta(y)\in \mbox{Ran}(\theta+\frac{a}{2})$ since $\theta$ is continuous, i.e., there exists an $w\in[0,1]$ such that \begin{equation}\label{eq7}\theta(x)+\theta(y)=\theta(w)+\frac{a}{2}.\end{equation} Let \begin{equation}\label{eq8} w_{1}=(\theta+\frac{a}{2})^{(-1)}(\theta(w)+\frac{a}{2}).\end{equation} Thus, from the continuity of $\theta$ we have \begin{equation}\label{eq08}\theta(w)=\theta(w_{1}).\end{equation} Therefore, for each $x,y,z\in[0,1]$,
		\begin{align*}
			T_{sub}(T_{sub}(x,y),z)=&T_{sub}((\theta+\frac{a}{2})^{(-1)}(\theta(x)+\theta(y)),z)\\=&(\theta+\frac{a}{2})^{(-1)}(\theta((\theta+\frac{a}{2})^{(-1)}(\theta(x)+\theta(y)))+\theta(z))\\
			=&(\theta+\frac{a}{2})^{(-1)}(\theta(w_{1})+\theta(z))\quad{\mbox{ by \eqref{eq7} and \eqref{eq8}}}\\
			=&(\theta+\frac{a}{2})^{(-1)}(\theta(w)+\theta(z))\quad{\mbox{ by \eqref{eq08}}}\\
			=&(\theta+\frac{a}{2})^{(-1)}(\theta(x)+\theta(y)+\theta(z)-\frac{a}{2}).
		\end{align*}
		Analogously, $T_{sub}(x,T_{sub}(y,z))=(\theta+\frac{a}{2})^{(-1)}(\theta(x)+\theta(y)+\theta(z)-\frac{a}{2})$ for each $x,y,z\in[0,1]$.
		Consequently, $T_{sub}(x,T_{sub}(y,z))=T_{sub}(T_{sub}(x,y),z)$ for each $x,y,z\in[0,1]$.
		\item[(\romannumeral3)] The monotonicity of $T_{sub}$. For $y,z\in[0,1]$ with $y\le z$, $\theta(y)\ge\theta(z)$ since $\theta$ is decreasing. This follows that  $T_{sub}(x,y)=(\theta+\frac{a}{2})^{(-1)}(\theta(x)+\theta(y))\le (\theta+\frac{a}{2})^{(-1)}(\theta(x)+\theta(z))=T_{sub}(x,z)$.
		\item[(\romannumeral4)] For each $x,y\in[0,1]$, we have that $T_{sub}(x,y)=(\theta+\frac{a}{2})^{(-1)}(\theta(x)+\theta(y))\le(\theta+\frac{a}{2})^{(-1)}(\theta(x)+\frac{a}{2})\le x$ since both $\theta$ and $\theta+\frac{a}{2}$ are decreasing. Similarly, $T_{sub}(x,y)\le y$. Therefore, $T_{sub}(x,y)\le\min\{x,y\}$.
(i), (ii), (iii) and (iv) mean that	$T_{sub}$ is a t-subnorm.	
\end{itemize}
		Step 2. We prove that $\mathcal{F}$ is a (pseudo) automorphism.

In fact, it is easy to see that $\mathcal{F}$ is continuous and increasing since both $\theta$ and $\vartheta$ are continuous and decreasing.
			From $\mathcal{F}(x)=\vartheta(\theta(x)+\frac{a}{2})$, we have that $\mathcal{F}(x)=0 \mbox{ if and only if } \vartheta(\theta(x)+\frac{a}{2})=0$. Moreover, by Theorem \ref{theorem3.2}, $\vartheta(\theta(x)+\frac{a}{2})=0\mbox{ if and only if } \theta(x)=\infty\mbox{ if and only if }\ x=0$ since $O_{\theta,\vartheta}$ is an overlap function. Then $\mathcal{F}(x)=0 \mbox{ if and only if }\ x=0$. Again from $\mathcal{F}(x)=\vartheta(\theta(x)+\frac{a}{2})$, we have that $\mathcal{F}(x)=1 \mbox{ if and only if } \vartheta(\theta(x)+\frac{a}{2})=1$. Moreover, from the hypothesis, it follows from Propositions \ref{proposition3.1} and \ref{proposition3.2} that we have $\vartheta(\theta(x)+\frac{a}{2})=1\mbox{ if and only if }\theta(x)+\frac{a}{2}\in[0,a]\mbox{ if and only if }\theta(x)=\frac{a}{2}\mbox{ if and only if }\ x=1$ since $O_{\theta,\vartheta}$ is an overlap function. Therefore, $\mathcal{F}(x)=1 \mbox{ if and only if }\ x=1$. Consequently, by Definition \ref{definition2.5}, $\mathcal{F}$ is a (pseudo) automorphism.\\
	Step 3. We prove that the overlap function $O_{\theta,\vartheta}$ can be obtained by a $(\mathcal{F},T_{sub})$-distortion for a (pseudo) automorphism $\mathcal{F}$ and a t-subnorm $T_{sub}$.	
		
		For any $x,y\in [0,1]$,
		\begin{align*}
			\mathcal{F}(T_{sub}(x,y))=&\vartheta(\theta(\theta+\frac{a}{2})^{(-1)}(\theta(x)+\theta(y))+\frac{a}{2})\quad{\mbox{ by \eqref{eq06} and \eqref{eq07}}}\\
			=&\vartheta(\theta(\theta+\frac{a}{2})^{(-1)}(\theta(w)+\frac{a}{2})+\frac{a}{2})\quad\mbox{by \eqref{eq7}}\\
			=&\vartheta(\theta(w_{1})+\frac{a}{2})\quad\mbox{by \eqref{eq8}}\\
			=&\vartheta(\theta(w)+\frac{a}{2})\quad{\mbox{ by \eqref{eq08}}}\\
			=&\vartheta(\theta(x)+\theta(y))\quad\mbox{by \eqref{eq7}}\\
			=&O_{\theta,\vartheta}(x,y),
		\end{align*}
i.e., $O_{\theta,\vartheta}(x,y)=\mathcal{F}(T_{sub}(x,y))$.
		\end{proof}

The following theorem gives two descriptions that an overlap function $O_{\theta,\vartheta}$ additively generated by the pair ($\theta$,$\vartheta$) can be obtained by a $(\mathcal{F},T)$-distortion for a (pseudo) automorphism $\mathcal{F}$ and a t-norm $T$.
\begin{theorem}\label{theorem5.1}
	For a function $O:[0,1]^{2}\rightarrow[0,1]$ the following are equivalent:
	\begin{itemize}
		\item[(1)] There exist a continuous and strictly decreasing function $\theta:[0,1]\rightarrow[0,\infty]$ and a continuous and decreasing function $\vartheta:[0,\infty]\rightarrow[0,1]$ such that the function $O$ is an overlap function additively generated by the pair ($\theta$,$\vartheta$) satisfying one of the following two conditions: for a given $a\in [0,\infty)$,
\begin{itemize}
		\item[(a)] $\theta(x)=\frac{a}{2}$ if and only if $x=1$;
		\item[(b)] $\vartheta(x)=1$ if and only if $x\in[0,a]$.
\end{itemize}
		\item[(2)] The function $O$ is an overlap function obtained by a ($\mathcal{F}$,$T$)-distortion for a (pseudo) automorphism $\mathcal{F}$ and a strict t-norm T;
		\item[(3)] There exist a strictly increasing bijection $\varphi:[0,1]\rightarrow[0,1]$ and a (pseudo) automorphism $\mathcal{H}$ such that $O(x,y)=\mathcal{H}(\varphi(x)\varphi(y))$ for all $(x,y)\in[0,1]^{2}$.
	\end{itemize}
\end{theorem}
\begin{proof}
$(1)\Rightarrow(2):$
	By Proposition \ref{proposition5.1} and its proof, the overlap function $O$ can be obtained by a $(\mathcal{F},T_{sub})$-distortion for a (pseudo) automorphism $\mathcal{F}$ and a t-subnorm $T_{sub}$ where $T_{sub}(x,y)=(\theta+\frac{a}{2})^{(-1)}(\theta(x)+\theta(y))$ for any $x,y\in [0,1]$. From the hypothesis, by Proposition \ref{proposition3.1} we always get that the condition (a) holds since $O$ is the overlap function. It is obvious that $T_{sub}$ is a continuous t-subnorm since $\theta$ is a continuous and strictly decreasing function, moreover,
	\begin{align*}
		T_{sub}(x,1)=&(\theta+\frac{a}{2})^{(-1)}(\theta(x)+\theta(1))\\
		=&(\theta+\frac{a}{2})^{(-1)}(\theta(x)+\frac{a}{2})\\
		=&x.
	\end{align*}
	Then $T_{sub}$ is a continuous t-norm. The rest of the proof are split into two parts as below.
	
	Part A. $T_{sub}$ is Archimedean.

 By Definition \ref{de2.2} we have
	\begin{align}\label{eq09}
		x_{T_{sub}}^{(n)}=&(\theta+\frac{a}{2})^{(-1)}(n\theta(x)-(n-2)\frac{a}{2})\nonumber\\
		=&\sup\{z\in[0,1]\mid\theta(z)+\frac{a}{2}>n\theta(x)-(n-2)\frac{a}{2}\}
	\end{align}for any $x\in(0,1)$.
	 Since $n\theta(x)-(n-2)\frac{a}{2}=n(\theta(x)-\frac{a}{2})+a$ and $\theta(x)>\frac{a}{2}$ for any $x\in(0,1)$, we get that $\lim\limits_{n\to\infty}n(\theta(x)-\frac{a}{2})+a=\infty$. This together with \eqref{eq09} yields $\lim\limits_{n\to\infty}x_{T_{sub}}^{(n)}=0$. Thus by Theorem 2.12 of \cite{EP2000}, $T_{sub}$ is Archimedean.

Part B. $T_{sub}$ is a t-norm with no zero divisors.

Indeed, by Theorem \ref{theorem3.2} the condition (1) of Theorem \ref{theorem3.1} holds since $O$ is an overlap function. Thus, we get $xy=0$ whenever $\sup\{z\in[0,1]\mid\theta(z)+\frac{a}{2}>\theta(x)+\theta(y)\}=0$.
Then it follows from
	 \begin{align*}
	 	T_{sub}(x,y)=&(\theta+\frac{a}{2})^{(-1)}(\theta(x)+\theta(y))\\
	 	=&\sup\{z\in[0,1]\mid\theta(z)+\frac{a}{2}>\theta(x)+\theta(y)\}
	 		 \end{align*} that, $T_{sub}(x,y)=0$ if and only $xy=0$. Thus $T_{sub}$ is a continuous t-norm with no zero divisors.

Therefore, $T_{sub}$ is a continuous Archimedean t-norm with no zero divisors, i.e., $T_{sub}$ is strict by Theorem 2.18 of \cite{EP2000}.
	
	$(2)\Rightarrow(1):$
	Suppose that $O(x,y)=\mathcal{F}(T(x,y))$ for a (pseudo) automorphism $\mathcal{F}$ and a strict t-norm $T$. Thus, from Theorem 5.1 of \cite{EP2000}, we know that $T$ has a continuous additive generator $t:[0,1]\rightarrow[0,\infty]$ where $t$ is a continuous and strictly decreasing function satisfying $t(1)=0$. Define $\theta:[0,1]\rightarrow[0,\infty]$ by \begin{equation}\label{equation14}
		\theta(x)=t(x)+\frac{a}{2}
	\end{equation} for any $x\in[0,1]$.	Then from Corollary 3.30 of \cite{EP2000}, $\theta$ is continuous and strictly decreasing and satisfies the following two conditions:
	\begin{itemize}
		\item[(i)] $\theta(x)=\frac{a}{2}\ \text{if and only if}\ x=1;$
		\item[(ii)] $\theta(x)=\infty\ \text{if and only if}\  x=0.$
	\end{itemize}
	Define $\vartheta:[0,\infty]\rightarrow[0,1]$ by
	\begin{equation}\label{equation(13)}
		\vartheta(x)=\mathcal{F}((\theta+\frac{a}{2})^{(-1)}(x)).
	\end{equation} for any $x\in[0,\infty]$. Then, $\vartheta$ is a continuous and decreasing function since $\mathcal{F}$ is increasing and continuous, and $(\theta+\frac{a}{2})^{(-1)}$ is decreasing and continuous.
	
	By Definition \ref{definition2.5}, $\mathcal{F}((\theta+\frac{a}{2})^{(-1)}(x))=0$ if and only if $(\theta+\frac{a}{2})^{(-1)}(x)=0$. Then from (\ref{equation(13)}), $\vartheta(x)=0$ if and only if $(\theta+\frac{a}{2})^{(-1)}(x)=0$. On the other hand, from (ii), we have $\sup\{z\in[0,1]\mid\theta(z)+\frac{a}{2}>x\}=0$ if and only if $x=\infty$. Thus $\vartheta(x)=0$ if and only if $x=\infty$.

From (\ref{equation(13)}) and Definition \ref{definition2.5}, it is obvious that $\vartheta(x)=1$ if and only if $(\theta+\frac{a}{2})^{(-1)}(x)=1$. From (i), we have $\sup\{z\in[0,1]\mid\theta(z)+\frac{a}{2}>x\}=1$ if and only if $x\in[0,a]$. Thus $\vartheta(x)=1$ if and only if $x\in[0,a]$.

Therefore,
	\begin{align*}
		\vartheta(\theta(x)+\theta(y))=&\mathcal{F}((\theta+\frac{a}{2})^{(-1)}(\theta(x)+\theta(y)))\quad\mbox{by \eqref{equation(13)}}\\
		=&\mathcal{F}(\sup\{z\in[0,1]\mid\theta(z)+\frac{a}{2}>\theta(x)+\theta(y)\})\quad\mbox{by Definition \ref{de2.2}}\\
		=&\mathcal{F}(\sup\{z\in[0,1]\mid t(z)+a>t(x)+t(y)+a\})\quad\mbox{by \eqref{equation14}}\\
		=&\mathcal{F}(\sup\{z\in[0,1]\mid t(z)>t(x)+t(y)\})\\
		=&\mathcal{F}(t^{(-1)}(t(x)+t(y)))\quad\mbox{by Definition \ref{de2.2}}\\
		=&\mathcal{F}(T(x,y))\\
		=&O(x,y).
	\end{align*}
	$(2)\Leftrightarrow(3)$ is immediately from \cite{FZ2022}.
\end{proof}

\begin{example}
	\emph{For a given $a\in[0,\infty)$, let $\theta:[0,1]\rightarrow[0,\infty]$ and $\vartheta:[0,\infty]\rightarrow[0,1]$ be defined, respectively, by
		$$\theta(x)=\begin{cases}
		\frac{a}{2}+ln\frac{2-x}{x} &\mbox{if}\ x\neq0,\\
		\infty &\mbox{if}\ x=0
	\end{cases}$$
	for any $x\in [0,1]$ and
	$$\vartheta(x)=
	\begin{cases}
		1       &\mbox{if}\ x\in[0,a),\\
		\frac{4}{e^{2(x-a)}+2e^{x-a}+1}&\mbox{if}\ x\in[a,\infty]
	\end{cases}$$
	for any $x\in [0,\infty]$. Then, it is obvious that $\theta$ is a continuous and strictly decreasing function and, $\vartheta$ is a continuous and decreasing function.
	By a simple computation, we get that $$O_{\theta,\vartheta}(x,y)=\vartheta(\theta(x)+\theta(y))=(\frac{xy}{2-x-y+xy})^{2}.$$
	Hence, from Theorem \ref{theorem5.1}, the function $O_{\theta,\vartheta}$ is an overlap function obtained by a ($\mathcal{F},T$)-distortion for a (pseudo) automorphism $\mathcal{F}$ and a strict t-norm $T$. In fact, the Hamacher t-norm $T(x,y)=\frac{xy}{2-x-y+xy}$ is a strict t-norm and let $\mathcal{F}(x)=x^{2}$. Then
	\begin{equation*}
		O_{\theta,\vartheta}(x,y)=(\frac{xy}{2-x-y+xy})^{2}=\mathcal{F}(T(x,y)).
	\end{equation*}
	Meanwhile, from Theorem \ref{theorem5.1}, there exist a strictly increasing bijection $\varphi:[0,1]\rightarrow[0,1]$ and a (pseudo) automorphism $\mathcal{H}$ such that $O(x,y)=\mathcal{H}(\varphi(x)\varphi(y))$ for each $(x,y)\in[0,1]^{2}$. Indeed, let $\varphi(x)=\frac{x}{2-x}$ and $\mathcal{H}(x)=(\frac{2x}{1+x})^{2}$. Obviously, $\varphi$ is a strictly increasing bijection and $\mathcal{H}$ is an automorphism. Moreover,
	\begin{equation*}
		\mathcal{H}(\varphi(x)\varphi(y))
		=(\frac{xy}{2-x-y+xy})^{2}
	\end{equation*}
	for all $(x,y)\in[0,1]^{2}$.}
\end{example}

As a conclusion of this section, we investigate when an overlap function additively generated by the pair $(\theta,\vartheta)$ cannot be obtained by a $(\mathcal{F},T)$-distortion for a positive and continuous t-norm $T$ and a (pseudo) automorphism $\mathcal{F}$. First, we introduce two necessary lemmas.

\begin{lemma}[\cite{FZ2022}]\label{lemma5.1}
	Let $T$ be a positive continuous t-norm. Then, one of the following three cases is valid for $T$:
	\begin{itemize}
		\item[(1)] $T(x,y)=min\{x,y\};$
		\item[(2)] $T(x,y)$ is a strict t-norm;
		\item[(3)] There exists a countable family $\{[a_{\alpha},e_{\alpha}],T_{\alpha}\}$ such that $T$ is an ordinal sum of this family and each $T_{\alpha}$ is a continuous Archimedean t-norm, and $T_{\alpha}$ is a strict t-norm when $a_{\alpha}=0$.
	\end{itemize}
\end{lemma}

\begin{lemma}[\cite{FZ2022}]\label{lemma5.2}
	Let $O:[0,1]^2\rightarrow[0,1]$ be an overlap function and $\mathcal{F}:[0,1]\rightarrow[0,1]$ be a (pseudo) automorphism with $\mathcal{F}(x)=O(x,1)$. Then the following are equivalent:
	\begin{itemize}
		\item[(1)] The overlap function $O$ can be determined by a $(\mathcal{F},T)$-distortion for the (pseudo) automorphism $\mathcal{F}$ and a positive continuous t-norm $T$.
		\item[(2)] One of the following three cases is valid for the overlap function $O$:
		\begin{itemize}
			\item[(\romannumeral1)] $O(x,y)=\min\{\mathcal{F}(x),\mathcal{F}(y)\}$;
			\item[(\romannumeral2)] There exists a continuous and strictly decreasing function $t:[0,1]\rightarrow[0,\infty]$ with $t(1)=0$ and $t(0)=\infty$ such that $$O(x,y)=\mathcal{F}(t^{(-1)}(t(x)+t(y)))$$ for all $(x,y)\in[0,1]^2$;
			\item[(\romannumeral3)] There exists a countable family $[a_{\alpha},e_{\alpha}]$ of non-overlapping, closed, non-trivial, proper subintervals of $[0,1]$ such that, for all $(x,y)\in[0,1]^{2}$, $$O(x,y)=\left\{
			\begin{aligned}
				&\mathcal{F}(h_{\alpha}^{-1}(h_{\alpha}(x)+h_{\alpha}(y)))       &if\ (x,y)\in[a_{\alpha},e_{\alpha}]^{2},\\
				&\min\{\mathcal{F}(x),\mathcal{F}(y)\}&otherwise,
			\end{aligned}	
			\right.$$
			where $h_{\alpha}:[a_{\alpha},e_{\alpha}]\rightarrow[0,\infty]$ is a continuous and strictly decreasing function which satisfies $h_{\alpha}(e_{\alpha})=0$ for each $\alpha$ and $h_{\alpha}(a_{\alpha})=\infty$ when $a_{\alpha}=0$.
		\end{itemize}
	\end{itemize}
\end{lemma}

Then we have the following proposition.
\begin{proposition}\label{proposition5.2}
	Let $\theta:[0,1]\rightarrow[0,\infty]$ be a continuous and decreasing but not strict function and $\vartheta:[0,\infty]\rightarrow[0,1]$ be a continuous and decreasing function. Let $O_{\theta,\vartheta}:[0,1]^2\rightarrow[0,1]$ be an overlap function additively generated by the pair $(\theta,\vartheta)$ that satisfies one of the following two conditions: for a given $a\in [0,\infty)$,
	\begin{itemize}
		\item[(1)] $\theta(x)=\frac{a}{2}$ if and only if $x=1$;
		\item[(2)] $\vartheta(x)=1$ if and only if $x\in[0,a]$.
	\end{itemize}
Then the overlap function $O_{\theta,\vartheta}$ cannot be obtained by a $(\mathcal{F},T)$-distortion for a positive and continuous t-norm $T$ and a (pseudo) automorphism $\mathcal{F}$.
\end{proposition}
\begin{proof}
	Suppose that the overlap function $O_{\theta,\vartheta}$ can be obtained by a $(\mathcal{F},T)$-distortion for a positive and continuous t-norm $T$ and a (pseudo) automorphism $\mathcal{F}$. Then by Lemma \ref{lemma5.1}, there are three cases as follows.
	
	Case 1. If $T(x,y)=\min\{x,y\}$, then $O_{\theta,\vartheta}(x,y)=\min\{\mathcal{F}(x),\mathcal{F}(y)\}$. In particular, \begin{equation}\label{eq015} O_{\theta,\vartheta}(x,x)=\vartheta(\theta(x)+\theta(x))=\mathcal{F}(x)\end{equation} for any $x\in(0,1)$. From the hypothesis, it follows from Proposition \ref{proposition3.1} that we have $\theta(x)=\frac{a}{2}$ if and only if $x=1$ since $O_{\theta,\vartheta}$ is an overlap function. Thus \begin{equation}\label{eq016}O_{\theta,\vartheta}(x,1)=\vartheta(\theta(x)+\frac{a}{2})=\mathcal{F}(x)\end{equation} and, $\theta(x)>\frac{a}{2}$ for any $x\in(0,1)$ since $\theta$ is decreasing. The inequality means that $2\theta(x)>\theta(x)+\frac{a}{2}$. Hence there exists a $y$ with $0<y<x$ such that $\theta(x)+\frac{a}{2}<\theta(y)+\frac{a}{2}<2\theta(x)$. This together with \eqref{eq015} and \eqref{eq016} results in $\mathcal{F}(x)=\vartheta(2\theta(x))=\vartheta(\theta(y)+\frac{a}{2})=\vartheta(2\theta(y))=\mathcal{F}(y)$. By induction, we can get that $\mathcal{F}(u)=\mathcal{F}(x)$ for any $u\in(0,x)$. Then $\mathcal{F}$ is a constant over $(0,x)$ and $\mathcal{F}(0)=0$, which violate the continuity of $\mathcal{F}$.
	
	Case 2. If $T(x,y)$ is a strict t-norm, then by Theorem \ref{theorem5.1}, $\theta$ is a strict function, a contradiction.
	
	Case 3. If there exists a countable family $\{[a_{\alpha},e_{\alpha}],T_{\alpha}\}$ such that $T$ is an ordinal sum of this family and each $T_{\alpha}$ is a continuous Archimedean t-norm, and $T_{\alpha}$ is a strict t-norm when $a_{\alpha}=0$, then from Lemma \ref{lemma5.2}, there exists a countable family $[a_{\alpha},e_{\alpha}]$ of non-overlapping, closed, non-trivial, proper subintervals of $[0,1]$ such that, for all $(x,y)\in[0,1]^{2}$,
	\begin{equation}\label{equation4}
		O(x,y)=\left\{
		\begin{aligned}
			&\mathcal{F}(h_{\alpha}^{(-1)}(h_{\alpha}(x)+h_{\alpha}(y)))       &\mbox{if}\ (x,y)\in[a_{\alpha},e_{\alpha}]^{2}\\
			&\min\{\mathcal{F}(x),\mathcal{F}(y)\}&\mbox{otherwise},
		\end{aligned}	
		\right.
	\end{equation}
	where $h_{\alpha}:[a_{\alpha},e_{\alpha}]\rightarrow[0,\infty]$ is a continuous and strictly decreasing function which satisfies $h_{\alpha}(e_{\alpha})=0$ for each $\alpha$ and $h_{\alpha}(a_{\alpha})=\infty$ when $a_{\alpha}=0$.
	
	It is easy to see that the interval $[0,1]$ can be represented as a union of a countable family of pairwise disjoint intervals $(J_{\alpha})_{\alpha\in A}$ where, for each $\alpha\in A$, either $J_{\alpha}=[u_{\alpha},v_{\alpha}]$ or $J_{\alpha}=(u_{\alpha},v_{\alpha})$ for suitable $u_{\alpha},v_{\alpha}\in[0,1]$. From Theorem \ref{theorem3.2}, we have $\theta(x)=\infty$ if and only if $x=0$ since $O_{\theta,\vartheta}$ is an overlap function and, from hypothesis and Proposition \ref{proposition3.1} we get $\theta(x)=\frac{a}{2}$ if and only if $x=1$. Then there exist suitable $u_{\alpha},v_{\alpha}$ such that $u_{\alpha}\ne0$ and $\theta(u_{\alpha})>\theta(v_{\alpha})$ and $[u_{\alpha},v_{\alpha}]$ or $(u_{\alpha},v_{\alpha})\in(J_{\alpha})_{\alpha\in A}$. Hence, for some fixed $x$ with $0<x<u_{\alpha}$, we have $\vartheta(\theta(x)+\theta(v_{\alpha}))=\vartheta(\theta(x)+\theta(u_{\alpha}))=\mathcal{F}(x)$ by \eqref{equation4}. Because $\theta$ is decreasing and continuous, there exists an $x_{0}$ with $0<x_{0}<x$ such that $\theta(x)+\theta(v_{\alpha})<\theta(x_{0})+\theta(v_{\alpha})<\theta(x)+\theta(u_{\alpha})$ and, $\mathcal{F}(x_{0})=\vartheta(\theta(x_{0})+\theta(v_{\alpha}))$ from \eqref{equation4}. Then, $\mathcal{F}(x_{0})=\vartheta(\theta(x_{0})+\theta(v_{\alpha}))=\vartheta(\theta(x+\theta(v_{\alpha})))=\mathcal{F}(x)$. Therefore, by induction, $\vartheta(\theta(u)+\theta(v_{\alpha}))=\mathcal{F}(u)=\mathcal{F}(x)\ne0$ for any $u\in(0,x)$, which violates the continuity of $\mathcal{F}$ since $\mathcal{F}(x)=0$ if and only if $x=0$.
	
	This completes the proof.		
\end{proof}

\section{Additive generator pairs of grouping functions}\label{section6}
In this section, we introduce the concept of an additive generator pair of grouping functions and some related results whose proofs are omitted since grouping functions are given as dual operations of overlap functions.

First, we give the definition of an additive generator pair of a grouping function.

\begin{definition}
	\emph{Let $t:[0,1]\rightarrow[0,\infty]$ and $s:[0,\infty]\rightarrow[0,1]$ be continuous and increasing functions, respectively. If the bivariate function $G_{t,s}:[0,1]^{2}\rightarrow[0,1]$, given by
		$$G_{t,s}(x,y)=s(t(x)+t(y)),$$
		is a grouping function, then $(t,s)$ is called an additive generator pair of the grouping function $G_{t,s}$ and $G_{t,s}$ is said to be additively generated by the pair $(t,s)$.}
\end{definition}

The following theorem characterizes the functions $t:[0,1]\rightarrow[0,\infty]$ and $s:[0,\infty]\rightarrow[0,1]$, respectively, whenever $G_{t,s}$ is a grouping function additively generated by the pair ($t$,$s$).
\begin{theorem}
	Let $t:[0,1]\rightarrow[0,\infty]$ and $s:[0,\infty]\rightarrow[0,1]$ be continuous and increasing functions, respectively. Let $G_{t,s}:[0,1]^2\rightarrow[0,1]$ be a grouping function additively generated by the pair $(t,s)$. Then the following statements hold:
	\begin{itemize}
		\item[(1)] $t(x)=\infty$ if and only if $x=1$;
		\item[(2)] $s(x)=1$ if and only if $x=\infty$.
	\end{itemize}
\end{theorem}

The following theorem shows the conditions that a pair $(t,s)$ can additively generate grouping functions.
\begin{theorem}\label{theorem6.1}
	For a given $a\in[0,\infty)$, let $t:[0,1]\rightarrow[0,\infty]$ and $s:[0,\infty]\rightarrow[0,1]$ be continuous and increasing functions such that
	\begin{itemize}
		\item[1.] $t(x)=\infty$ if and only if $x=1$;
		\item[2.] $t(x)=\frac{a}{2}$ if and only if $x=0$;
		\item[3.] $s(x)=0$ if and only if $x\in[0,a]$;
		\item[4.] $s(x)=1$ if and only if $x=\infty$.
	\end{itemize}
	Then, the function $G_{t,s}:[0,1]^{2}\rightarrow[0,1]$, defined by
	\begin{equation*}
		G_{t,s}(x,y)=s(t(x)+t(y)),
	\end{equation*}
	is a grouping function.
\end{theorem}

The following two propositions reveal that the conditions (2) and (3) of Theorem \ref{theorem6.1} are closely associated whenever $G_{t,s}:[0,1]^{2}\rightarrow[0,1]$ is a grouping function additively generated by the pair $(t,s)$.
\begin{proposition}
	For a given $a\in[0,\infty)$, let $t:[0,1]\rightarrow[0,\infty]$ and $s:[0,\infty]\rightarrow[0,1]$ be continuous and increasing functions satisfying
	\begin{itemize}
		\item[1.] $s(x)=0$ if and only if $x\in[0,a]$;
		\item[2.] a binary function $G_{t,s}:[0,1]^{2}\rightarrow[0,1]$ with $G_{t,s}(x,y)=s(t(x)+t(y))$ is a grouping function.
	\end{itemize}
	Then, $t(x)=\frac{a}{2}$ if and only if $x=0$.
\end{proposition}

\begin{proposition}
	For a given $a\in[0,\infty)$, let $t:[0,1]\rightarrow[0,\infty]$ and $s:[0,\infty]\rightarrow[0,1]$ be continuous and increasing functions satisfying
	\begin{itemize}
		\item[1.] $t(x)=\frac{a}{2}$ if and only if $x=0$;
		\item[2.] a binary function $G_{t,s}:[0,1]^{2}\rightarrow[0,1]$ with $G_{t,s}(x,y)=s(t(x)+t(y))$ is a grouping function.
	\end{itemize}
	Then, $x\in[0,a]$ if and only if $s(x)=0$.
\end{proposition}

In what follows, we shall investigate the conditions under which a grouping function additively generated by the pair ($t$,$s$) is a t-conorm.
\begin{proposition}
	Let $t:[0,1]\rightarrow[0,\infty]$ and $s:[0,\infty]\rightarrow[0,1]$ be continuous and increasing functions, respectively. Let $G_{t,s}:[0,1]^{2}\rightarrow[0,1]$ be a grouping function additively generated by the pair $(t,s)$ with $t(0)=\frac{a}{2}$ for some $a\in[0,\infty)$. Then the following statements are equivalent:
	\begin{itemize}
		\item[(1)] $s\circ(t+\frac{a}{2})=\emph{id}_{[0,1]}$;
		\item[(2)] $0$ is a neutral element of $G_{t,s}$.
	\end{itemize}
\end{proposition}

\begin{proposition}
	For a given $a\in[0,\infty)$, let $t:[0,1]\rightarrow[0,\infty]$ be a continuous and strictly increasing function such that:
	\begin{itemize}
		\item[(1)] $t(x)=\infty$ if and only if $x=1$;
		\item[(2)] $t(x)=\frac{a}{2}$ if and only if $x=0$.
	\end{itemize}
	If a function $s:[0,\infty]\rightarrow[0,1]$ satisfies $s(x)=(t+\frac{a}{2})^{(-1)}(x)$ for any $x\in [0,1]$, then a binary function $G_{t,s}:[0,1]^{2}\rightarrow[0,1]$ with $G_{t,s}(x,y)=s(t(x)+t(y))$ is a positive t-conorm.
\end{proposition}

\begin{proposition}
	Let $t:[0,1]\rightarrow[0,\infty]$ and $s:[0,\infty]\rightarrow[0,1]$ be two continuous and increasing functions, and $G_{t,s}:[0,1]^{2}\rightarrow[0,1]$ be a grouping function additively generated by the pair $(t,s)$ with $0$ as neutral element. If for a given $a\in [0,\infty)$, $t$ satisfies the condition (2) of Theorem \ref{theorem6.1}, then $G_{t,s}$ is associative.
\end{proposition}

\begin{theorem}
	Let $t:[0,1]\rightarrow[0,\infty]$ and $s:[0,\infty]\rightarrow[0,1]$ be two continuous and increasing functions and, $G_{t,s}:[0,1]^{2}\rightarrow[0,1]$ be a grouping function additively generated by the pair $(t,s)$. If for a given $a\in [0,\infty)$, $t$ satisfies the condition (2) of Theorem \ref{theorem6.1}, then the following three statements are equivalent:
	\begin{itemize}
		\item[(1)] $G_{t,s}$ is a t-conorm;
		\item[(2)] $0$ is a neutral element of $G_{t,s}$;
		\item[(3)] $s\circ(t+\frac{a}{2})=\emph{id}_{[0,1]}$.
	\end{itemize}
\end{theorem}

In the following, we give the conditions that a grouping function additively generated by the pair ($t$,$s$) can be obtained by a distortion of a t-conorm and a (pseudo) automorphism.
\begin{proposition}
	Let $\mathcal{F}:[0,1]\rightarrow[0,1]$ be a pseudo automorphism. Then for every positive and continuous t-conorm $S:[0,1]^{2}\rightarrow[0,1]$, the function $G_{\mathcal{F},S}:[0,1]^{2}\rightarrow[0,1]$, given by
	\begin{equation*}
		G_{\mathcal{F},S}(x,y)=\mathcal{F}(S(x,y)),
	\end{equation*}
	is a grouping function.
	
	Where $G_{\mathcal{F},S}$ is called a grouping function obtained by the distortion of the t-conorm $S$ by the pseudo automorphism $\mathcal{F}$ or a grouping function obtained by a $(\mathcal{F},S)$-distortion.
\end{proposition}

\begin{proposition}
	Let $t:[0,1]\rightarrow[0,\infty]$ and $s:[0,\infty]\rightarrow[0,1]$ be two continuous and increasing functions and, $G_{t,s}:[0,1]^{2}\rightarrow[0,1]$ be a grouping function additively generated by the pair $(t,s)$ that satisfies one of the following two conditions: for a given $a\in [0,\infty)$,
	\begin{itemize}
		\item[(1)] $t(x)=\frac{a}{2}$ if and only if $x=0$;
		\item[(2)] $s(x)=0$ if and only if $x\in[0,a]$.
	\end{itemize}
	Then the grouping function $G_{t,s}$ can be obtained by a $(\mathcal{F},T_{super})$-distortion for a (pseudo) automorphism $\mathcal{F}$ and a t-superconorm $T_{super}$.
\end{proposition}

The next theorem gives two descriptions that a grouping function $G_{t,s}$ additively generated by the pair $(t,s)$ can be obtained by a $(\mathcal{F},S)$-distortion for a (pseudo) automorphism $\mathcal{F}$ and a t-conorm $S$.
\begin{theorem}
	For a function $G:[0,1]^{2}\rightarrow[0,1]$ the following are equivalent:
	\begin{itemize}
		\item[(1)] There exist a continuous and strictly increasing function $t:[0,1]\rightarrow[0,\infty]$ and a continuous and increasing function $s:[0,\infty]\rightarrow[0,1]$ such that the function $G$ is a grouping function additively generated by the pair ($t$,$s$) satisfying one of the following two conditions: for a given $a\in [0,\infty)$,
		\begin{itemize}
			\item[(a)] $t(x)=\frac{a}{2}$ if and only if $x=0$;
			\item[(b)] $s(x)=0$ if and only if $x\in[0,a]$.
		\end{itemize}
		\item[(2)] The function $G$ is a grouping function obtained by a ($\mathcal{F}$,$S$)-distortion for a (pseudo) automorphism $\mathcal{F}$ and a strict t-conorm S;
		\item[(3)] There exist a strictly increasing bijection $\varphi:[0,1]\rightarrow[0,1]$ and a (pseudo) automorphism $\mathcal{H}$ such that $G(x,y)=\mathcal{H}(\varphi(x)+\varphi(y)-\varphi(x)\varphi(y))$ for all $(x,y)\in[0,1]^{2}$.
	\end{itemize}
\end{theorem}

Finally, we investigate when a grouping function additively generated by the pair $(t,s)$ cannot be obtained by a $(\mathcal{F},S)$-distortion for a positive and continuous t-conorm $S$ and a (pseudo) automorphism $\mathcal{F}$.
\begin{proposition}
	Let $t:[0,1]\rightarrow[0,\infty]$ be a continuous and increasing but not strict function and $s:[0,\infty]\rightarrow[0,1]$ be a continuous and increasing function. Let $G_{t,s}:[0,1]^2\rightarrow[0,1]$ be a grouping function additively generated by the pair $(t,s)$ that satisfies one of the following two conditions: for a given $a\in [0,\infty)$,
	\begin{itemize}
		\item[(1)] $t(x)=\frac{a}{2}$ if and only if $x=0$;
		\item[(2)] $s(x)=0$ if and only if $x\in[0,a]$.
	\end{itemize}
	Then the grouping function $G_{t,s}$ cannot be obtained by a $(\mathcal{F},S)$-distortion for a positive and continuous t-conorm $S$ and a (pseudo) automorphism $\mathcal{F}$.
\end{proposition}

\section{Conclusions}\label{section7}
The article mainly explored the relationship between functions $\theta$ and $\vartheta$ in an overlap function additively generated by a relaxing boundary additive generator pair ($\theta$,$\vartheta$) and characterized the overlap function additively generated by the additive generator pair ($\theta$,$\vartheta$). Compared with the known definition, see \cite{GP2016}, of an additive generator pair for overlap functions, our functions $\theta$ and $\vartheta$ in the additive generator pair of an overlap function are more general with a relaxed boundary conditions while the property of overlap functions can be preserved. Clearly, multiplicative generator pairs of overlap and grouping functions can be discussed analogously. In the future, it is interesting to consider whether the continuity of functions $\theta$ or $\vartheta$ in an overlap function additively generated by an additive generator pair ($\theta$,$\vartheta$) can be relaxed or not.


\section*{Declaration of competing interest}
The authors declare that they have no known competing financial interests or personal relationships that could have appeared to influence the work reported in this article.

\end{document}